\documentclass{tip}
\newcommand{\seqnum}[1]{\href{https://oeis.org/#1}{\rm \underline{#1}}}
\theoremstyle{plain}
\newtheorem{theorem}{Theorem}

\newtheorem{lemma}[theorem]{Lemma}

\theoremstyle{definition}
\newtheorem{definition}[theorem]{Definition}

\theoremstyle{remark}

\newcolumntype{L}{>{$}l<{$}} 
\newcolumntype{C}{>{$}c<{$}} 
\newcommand{\autora}{R.~Zöllner}
\newcommand{\autorb}{K.~Handrich}
\title{\vspace{-2.5cm}\hrulefill\\\color{darkblue}\sffamily\LARGE\bfseries On the Number of Posets\\}
\author{\sffamily\bfseries \autora\orcidlink{0000-0002-3544-6622}${}^{a,b}$~\&~\autorb${}^{a,c}$}
\affil{\sffamily${}^a$Institute of Material Handling and Industrial Engineering, TU~Dresden,\newline01062 Dresden, Germany}
\affil{\sffamily${}^b$rico.zoellner@tu-dresden.de}
\affil{\sffamily${}^c$konrad.handrich@tu-dresden.de}
\date{\sffamily\today\\\hrulefill}

\begin{document}
	\maketitle
	\thispagestyle{plain}
	\vspace{-2.0em}
	
	\parbox[t]{0.4\linewidth}{\textbf{Keywords}
		\begin{itemize}
			\item partially ordered sets
			\item Hasse diagrams
			\item enumeraitve combinatorics
			\end{itemize}}
	\hfill
	\parbox[t]{0.5\linewidth}{\textbf{Abstract}\\	This paper presents combinatorial facts dealing with the number of unlabeled partially ordered sets (posets) refined by the number of arcs in the Hasse diagram (sequence \seqnum{A342447} in OEIS). The main result is that the differences with respect to the number of points in this sequence become stationary if the number of points is sufficiently high. These differences are proposed as the new sequence \seqnum{A376894}. In addition, the underlying combinatorial and graph theoretical arguments were used to extend some further OEIS sequences.}
	
	\vspace{2.0em} 
	\noindent\rule{\textwidth}{0.4pt}
	\vspace{-3.0em}
	
	\section{Introdction}\label{Intro}
	Partial orders undoubtedly belong to the fundamental structures of mathematics. Their relevance can hardly be overestimated. As a natural consequence, partially ordered sets -- briefly called posets -- have become the subject of numerous studies. One part of these studies is dedicated to combinatorial questions. Here, Stanley's two-volume encyclopedia \cite{Stanley1, Stanley2} is a never ending source of knowledge, especially chapter 3 of the first volume~\cite{Stanley1}. Although posets seem to be quite simple objects from the structural point of view, not even a closed formula is known to compute the number of posets for a given number of points\footnote{However, the asymptotic behavior is known due to Kleitman~\cite{Kleitman}. As a consequence, Stanley~\cite{Stanley1} provides for the number of isomorphism classes of posets $H(p)$ (i.e., unlabeled posets) the asymptotic formula: 
	\begin{equation}
		H(p)\ \sim 2^{\frac{p^2}{4}}~\text{for}~p \to \infty.
\end{equation}}.
Brinkmann et al.\ \cite{McKay} calculated these numbers up to $16$ points (see on OEIS the sequences \seqnum{A000112} for unlabeled posets, and \seqnum{A000798} for labeled posets). Regarding the OEIS database, there are several sequences dealing with unlabeled posets summarized in the following definition.
\begin{definition}\label{X13}
	For $p>0$ unlabeled points (i.e., up to isomorphism) and $a \geq 0$ arcs in the related Hasse diagram (see Definition~\ref{Definition8}) we define
	\begin{enumerate}[label=(\roman*)]
		\item $H(p)$: The number of posets (\seqnum{A000112}).
		\item $H(p,a)$: The number of posets ($H(p)$ refined by the number of arcs, \seqnum{A342447}).
		\item $H_c(p,a)$: The number of connected posets (\seqnum{A022017}).
		\item $H_0(p,a)$: The number of posets without isolated points.
		\item $\overset{\infty}{H}(a)$: The total number of unlabeled posets without isolated points and $a$ arcs in the Hasse diagram, i.e. $\overset{\infty}{H}$ are the convergent down rows of $H(p,a)$ (\seqnum{A022016}).
	\end{enumerate}
\end{definition}
	The aim of the present paper is to work out some features of \seqnum{A342447}, where we used the ``nauty'' library \cite{Naughty} to compute the data shown in this paper. The main point is that the differences $H(p+1,a)-H(p,a)$ converge for sufficiently high $a$ and $p \geq \frac{3}{2}a$. We proposed these limit differences as a new sequence \seqnum{A376894} in OEIS. Some combinatorial arguments are employed to extend the OEIS sequence \seqnum{A022016} and parts of \seqnum{A342447}.
	
	Our paper is organized as follows: Section~\ref{Over} gives an overview of $H(p,a)$ and some regions in the $p$-$a$ plane, which is to be understood as a collection of basic facts. Hereby, the focus lies on connectivity. Section~\ref{Asym} contains the proof of the convergence of the above mentioned differences in $H(p,a)$. We summarize in Section~\ref{Summ}. Notational conventions and value tables are in the Appendices \ref{notations} and \ref{values}, respectively.
\begin{figure}[ht!]
	\centering
	\includegraphics[width=0.7\textwidth]{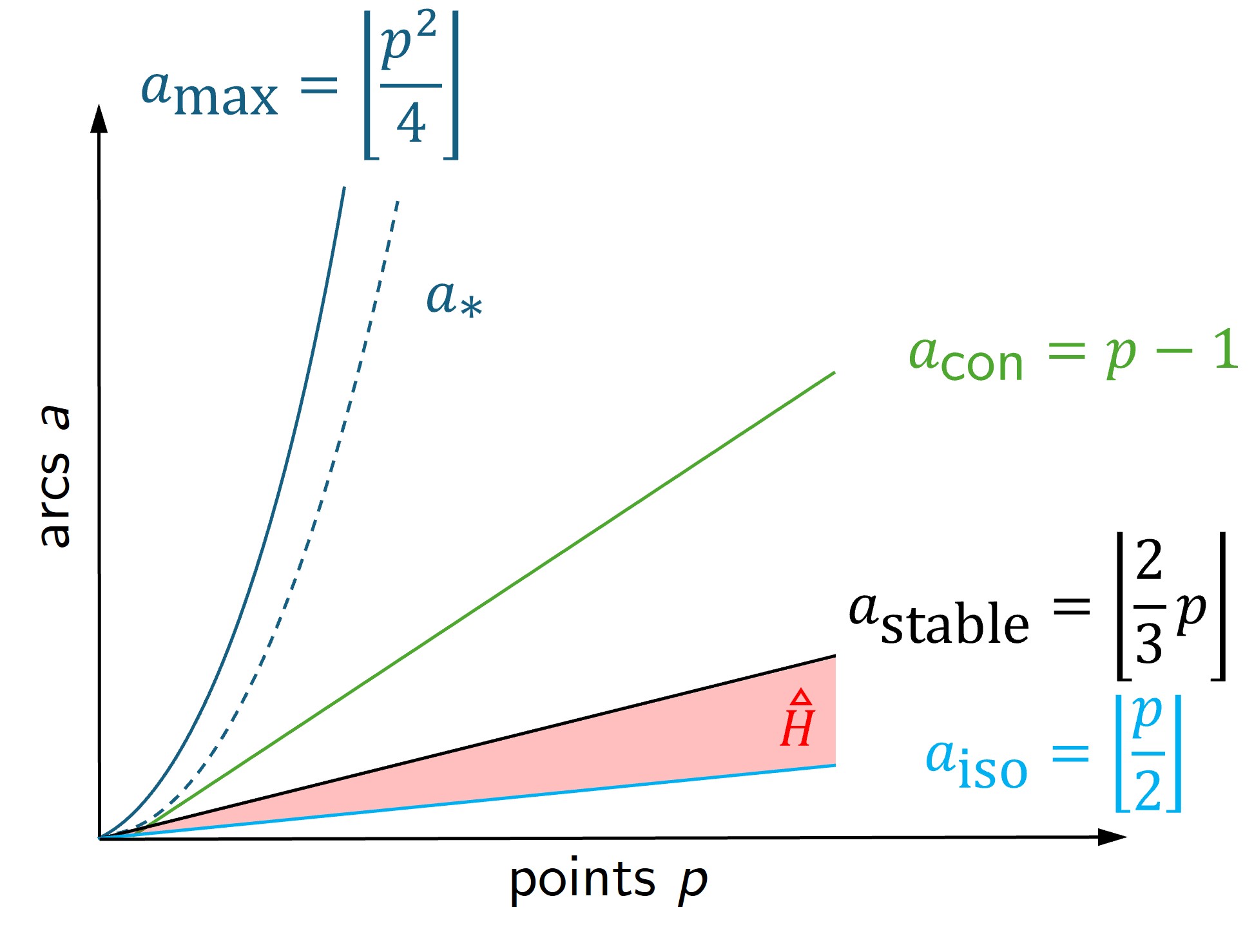}
	\caption{Schematic illustration of the relevant regions of the $p$-$a$~plane. For $a > a_{\max}$, there is no poset; for $a_{\max} (p) \geq a > a_\ast(p) = a_{\max}(p-1)$ all posets are connected, and below the $a_{\mathrm{con}}$ line all posets are disconnected. In addition, $H_0(p,a)$ is stationary in the cone spanned by the $a_{\mathrm{stable}}$ and $a_{\mathrm{iso}}$ line (leading to the $\overset{\Delta}{H}$ sequence \seqnum{A376894}); for $p \to \infty$ the number of posets converge to $\overset{\infty}{H}$.}
	\label{Fig1}
\end{figure}
\section{Overview}\label{Over}
In this section, we explore the landscape of posets regarding the number of points and the number of arcs in the Hasse diagram. Figure \ref{Fig1} displays the distinct regions and their borderlines. Besides the numbers $H(p,a)$ we pay attention to the number of posets without isolated points $H_0(p,a)$ and the number of connected posets $H_c(p,a)$, see Tables~\ref{normal} and \ref{conn} in Appendix~\ref{values}.

\begin{figure}[h]
	\centering
	\includegraphics[width=0.6\textwidth]{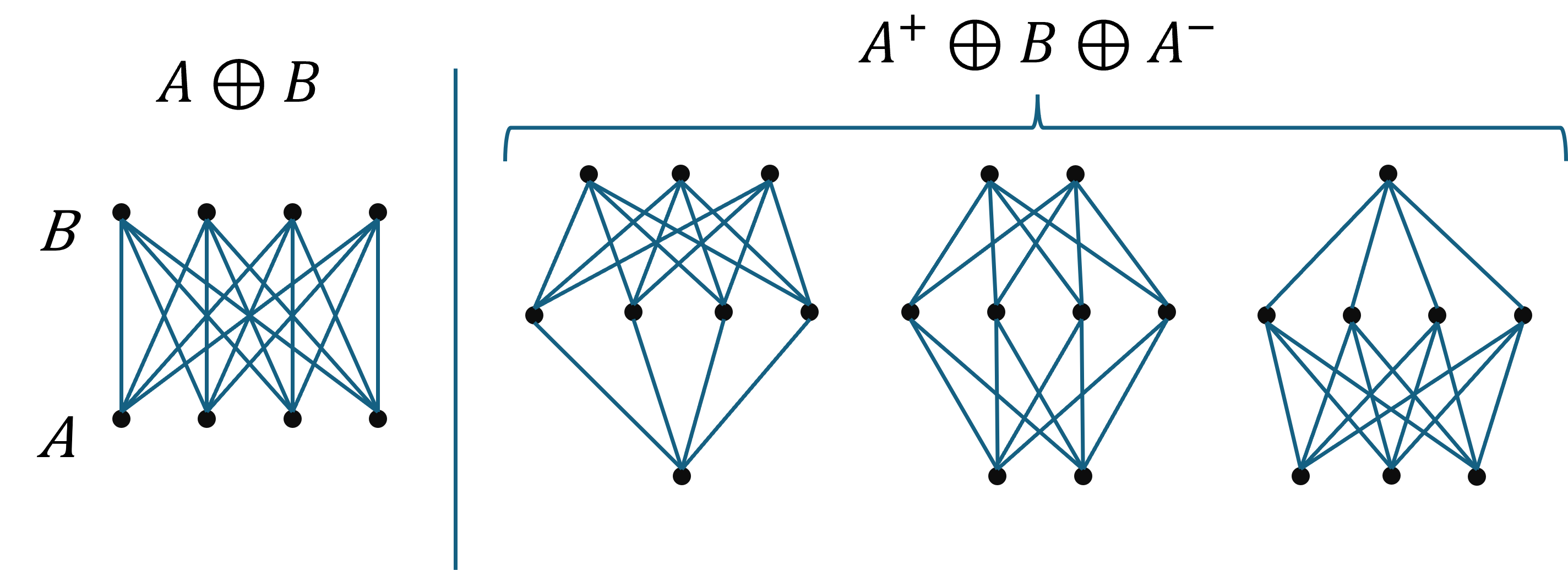}
	\caption{Illustration of the proof argument by the example with $p=8$, and $a_{\max}=16$, where $H(8,16)=4$. \textit{Left part:} This is basically the Tur\'an graph $T_2(8)$. \textit{Right part:} As the number of points from $A^-$ runs from $1$ to $3=\lfloor \frac{p}{2} \rfloor - 1$, three posets are generated from $T_2(8)$.}
	\label{X20}
\end{figure}

For a fixed number of points $p$, the situation is as follows. If the number of arcs $a$ is smaller than $p-1$, i.e., $a < p-1$ (below the $a_{\mathrm{con}}$ line in Figure~\ref{Fig1}), all posets are disconnected. If otherwise $a \geq p-1$, the total number $H(p,a)$ contains both connected and disconnected posets. Further increasing the arc number, we reach the dashed $a_\ast$ curve (see Figure~\ref{Fig1}) which is $a_{\max}(p-1)$ for $p>0$ and $a_\ast(0)=0$. Between $a_{\max}(p-1)$ and $a_{\max}(p)$, all posets are connected. Eventually above the $a_{\max}$ curve, there would be too many arcs to build a valid Hasse diagram, and thus, no poset can be found. The following lemma summarizes the description in a formal way:
\begin{lemma}\label{X4}
	\leavevmode\\\vspace{-5mm}
	\begin{enumerate}
		\item[(i)] If $a < p-1$, then $H_c(p,a)=0$.
		\item[(ii)] If $a > a_{\max}(p) = \lfloor\frac{p^2}{4}\rfloor$ then $H(p,a) = 0$. In addition, one has 
		\begin{equation}
			H(p,a_{\max}(p)) = \begin{cases}
				p, & \text{if}~p~\text{is odd}; \\
				\frac{p}{2}, & \text{if}~p~\text{is even}.
			\end{cases}
		\end{equation}
		\item[(iii)] If $a_{\max}(p-1) < a \leq a_{\max}(p)$, then $H(p,a) = H_c(p,a)$.
	\end{enumerate}
\end{lemma}
\begin{proof}
	Since all of those points are more or less known from graph theory, we only sketch the proof ideas. 
	\begin{enumerate}[label = (\roman*)]
		\item See the standard textbooks \cite[pp.\ 10]{Diestel}, \cite[pp.\ 73]{Bollobas}. 
		\item Due to Tur\'an's theorem \cite{Turan} (triangle-free case, also known as Mantel's theorem \cite{Mantel}), we have the bound for the number of edges as $\lfloor\frac{p^2}{4}\rfloor$. Furthermore, the only (undirected) graph reaching this bound is the Tur\'an graph $T_2(p)$, see the standard textbooks~\cite{Bollobas,Diestel,Aigner}. It remains to count all digraphs (Hasse diagrams), which are isomorphic to $T_2(p)$ after erasing the direction of the edges. To this end, note that $T_2(p)$ is bipartite with $T_2(p) = A \cup B$, where the two parts $A$ and $B$ have $\lfloor\frac{p}{2}\rfloor$ and $\lceil\frac{p}{2}\rceil$ points, respectively. If $p$ is even, we have $A \oplus B \cong B \oplus A$, yielding one poset. Note that the only way to orient $T_2(p)$ is to split either $A$ or $B$ into two parts. The splitting of one part (say $A$) into the two parts $A^+$, $A^-$ gives $\frac{p}{2}-1$ posets of the form $A^+ \oplus B \oplus A^-$ by choosing $A^+$ to have $1 \ldots\frac{p}{2}-1$ points (Figure~\ref{X20}). Thus, $1+\frac{p}{2} -1 = \frac{p}{2}$~posets if $p$ is even. Else if $p$ is odd, $A \oplus B$ and $ B \oplus A$ are two non-isomorphic posets. The analogous decompositions $A^+ \oplus B \oplus A^-$ and $B^+ \oplus A \oplus B^-$ lead to $\lfloor\frac{p}{2}\rfloor -1$ and $\lceil\frac{p}{2}\rceil -1$ combinations, respectively. In total, we have $2+\lfloor\frac{p}{2}\rfloor-1 + \lceil\frac{p}{2}\rceil-1 =p$ posets.
		\item Indirect: If the poset is disconnected, there are at least two non-empty connected components with $n$ and $p-n$ elements ($0<n<p$). Due to (ii) we conclude $a_{\max}(n) + a_{\max}(p-n) < a_{\max}(p)$ which is a contradiction.
	\end{enumerate}
\end{proof}
\noindent Now we turn to the disconnected region below the $a_{\mathrm{con}}$ line.
\begin{lemma}\label{X5}
	If $p>2a$, then $H_0(p,a)=0$.
\end{lemma}
\noindent Lemma~\ref{X5} justifies the Definition~\ref{X13}~(v) (see \seqnum{A022016}). The aim is to extend the sequence \seqnum{A022016} (as well as the disconnected part of \seqnum{A342447} by employing Lemma~\ref{X4}~(i) and (iii)). To this end, we consider decompositions of disconnected posets into distinct connected components without regard to isolated points.
\begin{definition}\label{X11}
	A \textit{distinct connected partition} of $(p,a)$ is a finite sequence $({p}_i, {a}_i, {f}_i)_{i=1}^{{s}}$ satisfying 
	\begin{enumerate}[label = (\roman*)]
		\item $\sum_{i=1}^{{s}}  {f}_i {p}_i = p$ and $\sum_{i=1}^{{s}}  {f}_i {a}_i = a$
		\item $a_i>0$ and $f_i>0$ for all $i=1 \ldots s$.
		\item If $i>j$, then either $p_i<p_j$ or $p_i=p_j$ with $a_i \leq a_j$ for all $i,j = 1 \ldots s$.
		\item All $(p_i,a_i)$ are pairwise distinct.
	\end{enumerate}
\end{definition}
\noindent Now the number $H_0(p,a)$ can be calculated using the well known formula for combinations with repetitions.
\begin{lemma}\label{X9}
	\begin{equation}
		H_0(p,a) = \sum_{(p_i, a_i, f_i)_{i=1}^{s}}^{}
		\prod_{i=1}^{s} 
		\binom{H_c({p}_i,{a}_i)+{f}_i -1}{H_c({p}_i,{a}_i)-1}
	\end{equation}
	where the sum is to be taken over all distinct connected partitions of $(p,a)$.
\end{lemma}
	\begin{proof}
	The number of combinations for choosing $m$ objects out of $n$ distinguishable objects allowing repetitions is equal to \begin{equation}
		\binom{n+m-1}{n-1}.
	\end{equation} 
	Regarding one element of a distinct connected partition, the combinations are equal to 
	\begin{equation} \binom{H_c(p_i,a_i)+f_i -1}{H_c(p_i,a_i)-1}. \end{equation} Since all elements are distinct, the combinations of one such a partition are simply given by the product of the binomial coefficients. And finally, all combinations have to be summed up.
\end{proof}
As an example, let us consider the calculation of $H_0(14,9)$. Note that there are six distinct connected partitions of $(14,9)$, which are listed in the following table:
\begin{table}[ht!]
	\centering
	\begin{tabular}[\textwidth]{r|rr|rrr|rr|rrr|rr|rr}
		\#& \multicolumn{2}{c|}{1}& \multicolumn{3}{c|}{2}&\multicolumn{2}{c|}{3}&\multicolumn{3}{c|}{4}&\multicolumn{2}{c|}{5}&\multicolumn{2}{c}{6} \\\midrule[1.5pt]
		$s$ &2&&3&&&2&&3&&&2&&2& \\
		$i$ &1&2&1&2&3&1&2&1&2&3&1&2&1&2\\ \midrule
		$p_i$&6&2&5&3&2&4&2&4&3&2&4&2&3&2\\
		$a_i$&5&1&4&2&1&3&1&3&2&1&4&1&2&1 \\
		$f_i$&1&4&1&1&3&2&3&1&2&2&1&5&4&1\\ \midrule
		$H_0(p_i,a_i)$&91&1&27&3&1&8&1&8&3&1&2&1&3&1\\
	\end{tabular}
\end{table}
The six summands according to Lemma~\ref{X9} are 
\begin{align}
	\binom{H_c(6,5)}{H_c(6,5) - 1}\cdot\binom{H_c(2,1)+3}{H_c(2,1) - 1}&=\binom{91}{90}\cdot\binom{4}{0}= 91 \\
	\binom{H_c(5,4)}{H_c(5,4) - 1}
	\cdot \binom{H_c(3,2)}{H_c(3,2) - 1}
	\cdot \binom {H_c(2,1)+2}{H_c(2,1) - 1}&= \binom{27}{26}\cdot\binom{3}{2}\cdot\binom{3}{0} = 81
\end{align}
and analogously $36$, $48$, $2$, and $15$. Hence $H_0(14,9)=91+81+36+48+2+15 = 273$ which we will later identify as $\overset{\Delta}{H}(5)$.
In summary, Lemma~\ref{X9} provides a method to compute the $H_0(p,a)$ numbers from some $H_c(\tilde{p},\tilde{a})$ numbers, where $\tilde{p} < p$ and $\tilde{a} < a$. This is useful to avoid complete enumerations. As shown in the next section, the asymptotic differences become stationary and can be computed up to a relatively high number of points.
\section{Asymptotic Differences}\label{Asym}
We observe that the differences $H(p,a)-H(p-1,a)$ in Table~\ref{normal} (which are the $H_0(p,a)$ numbers) stabilize between the $a_{\mathrm{iso}}$ line and the $a_{\mathrm{stable}}$ line. More precisely, in this region the $H_0(p,a)$ numbers are uniquely determined by the difference $p-2a$. In this section, we will prove that this is indeed the fact.
\begin{lemma}\label{X12}
	Let $n \geq 0$ be an arbitrary natural number. Then the expression $H_0(2a-n,a)$ does not depend on $a$ for $a \geq 2n$.
\end{lemma}
\begin{proof}
	Step 1: For $k>0$ and tuples of positive real numbers $(p_i)_{i=1}^k$, $(a_i)_{i=1}^k$ one has
	\begin{equation}\label{eq9}
		\min_{i=1\ldots k} \frac{p_i}{a_i} \leq \frac{p_1+\ldots+p_k}{a_1+\ldots+a_k} \leq \max_{i=1 \ldots k} \frac{p_i}{a_i}
	\end{equation}
	where equality holds if and only if $\frac{p_i}{a_i}$ are identical. Sometimes this inequality is referenced as inequality of the weighted arithmetic mean~\cite[p.\ 95]{Heuser}. 
	
\noindent	Step 2: From inequality~(\ref{eq9}), it follows that for $ p>\frac{3}{2}a$, each distinct connected partition of $(p,a)$ contains at least one pair\footnote{A pair is a (sub)poset with two points connected with an arc, and the two points have no further arcs.}. Indeed, if $p > \frac{3}{2}a$ then $\frac{3}{2} < \frac{p}{a} \leq \max_{i=1 \ldots k} \frac{p_i}{a_i}$. Thus, there is at least one $(p_i,a_i)$ with $ \frac{p_i}{a_i} > \frac{3}{2}$, i.e., a pair due to Lemma~\ref{X4}~(i). 
	
\noindent	Step 3: For the case $p=\frac{3}{2}a$ (and $a$ even) there is only one distinct connected partition without a pair, namely the one with $s=1$, $p_1=3$, $a_1=2$, and $f_1 = \frac{a}{2}$. Indeed, if $\frac{3}{2} \leq \max_{i=1\ldots k} \frac{p_i}{a_i}$ we have that $\frac{p_i}{a_i}=\frac{3}{2}$ for all $i$ if and only if equality holds. Again by Lemma~\ref{X4}~(i), we arrive at $p_i=3$ and $a_i=2$ for all $i$.
	
\noindent	Step 4: Consequently, all distinct connected partitions of $(p+2,a+1)$ can be obtained by adding one isolated pair to all distinct connected partitions of $(p,a)$.  Finally, $H_0(p+2,a+1) = H_0(p,a)$ for $p \geq \frac{3}{2}a$. Setting $p=2a-n$ completes the proof.
\end{proof}
\noindent	Due to Lemma~\ref{X12} the following definition makes sense.
\begin{definition}
	For $n \geq 0$ and $a \geq 2n$ set
	
	\begin{equation}\overset{\Delta}{H}(n) = H_0(2a-n,a) = H(2a-n,a) - H(2a-n-1,a)\end{equation}
\end{definition}
\noindent Going along the $a_{\mathrm{stable}}$ line, we evaluate the $\overset{\Delta}{H}$ sequence as $\overset{\Delta}{H}(n) = H_0(3n,2n)$, concretely

\begin{equation}\overset{\Delta}{H}(n) = 1,3,14,61,\ldots\end{equation}
See Table~\ref{Hdelta} for all elements of $\overset{\Delta}{H}$ computed so far (see as well the associated sequence in OEIS \seqnum{A376894}). The sequences $H$, $\overset{\infty}{H}$, and $\overset{\Delta}{H}$ are therefore related by

\begin{equation}H(2a-n-1,a) = \overset{\infty}{H}(a) - \sum_{i=0}^{n} \overset{\Delta}{H}(i)\end{equation}
for $n \geq 0$ and $a \geq 2n$ which allows the extension of \seqnum{A342447} and \seqnum{A022016} (Table~\ref{Hinf}).
	\section{Summary and Outlook}\label{Summ}
Besides the collection of some facts about the OEIS sequence \seqnum{A342447}, we have proved that the differences in \seqnum{A342447} become stationary in a triangle-shaped region of the plane spanned by the point axis and the arc axis. These stationary asymptotic differences are captured by the new OEIS sequence \seqnum{A376894}. In addition, we employed combinatorial facts to extend the OEIS sequences \seqnum{A022016} and the disconnected part of  \seqnum{A342447}, where we exploited the data gained with Brinkmann and McKay's poset generator.

Nevertheless, the number of arcs in the Hasse diagram and the connectivity are just two minor aspects of the rich structure of poset classes. For instance, let us take a glimpse at self-dual posets, i.e., posets that are isomorphism-invariant under reversion of the ordering (see e.g., sequences \seqnum{A133983} and \seqnum{A376633} (Table~\ref{selfdual}) in the OEIS database). Since a disconnected poset is self-dual if and only if all components are self-dual, one can also establish a sequence analogous to the  $\overset{\Delta}{H}$ numbers, say  $\overset{\Delta}{H}_{\mathrm{sd}}$, and by the same reasoning as above we obtain
$\overset{\Delta}{H}_{\mathrm{sd}} = 1,1,4,7,25,\ldots$

All in all, posets are still interesting research objects as they connect aspects of combinatorics and graph theory with many applications in practice.
	\section{Acknowledgments}
The work of Rico Z\"ollner is supported by DFG project $418727532$. The authors thank Jens-Uwe Grabowski for his helpful comments.

\begin{appendices}
	\section{Definitions and Notational Conventions} \label{notations}
	Let us recall some basic definitions and set up the notation (consistently with Stanley, chapter~3~\cite{Stanley1}).
	\begin{definition} \label{Definition8}
		A \textit{partially ordered set} (poset) is a set $P$ together with a binary relation $\leq$ fulfilling the following three conditions:
		\begin{enumerate}[label = (\roman*)]
			\item $x \leq x$
			\item $x \leq y $ and $y \leq x$ implies $y=x$
			\item $x \leq y$ and $y \leq z$ implies $x \leq z$ 
		\end{enumerate}
		for all $x,y,z \in P$. In addition, $x$ is said to be \textit{covered} by $y$ if $x \neq y$, $x \leq y$, and  there is no $w \in P$ with  $w \neq x $, $w \neq y $, and $x \leq w \leq y$. 
	\end{definition}
\noindent	The Hasse diagram is the digraph encoding the cover relation within a poset. Next, we define the two operations for posets.
	\begin{definition}
		Let $P,Q$ be disjoint posets.
		\begin{enumerate}[label = (\roman*)]
			\item The \textit{direct sum} $P+Q$ is defined as the poset on $P \cup Q$, where $x \leq y$ in $P + Q$ if either
			\begin{enumerate}
				\item[(a)] $x,y \in P$ and $x \leq y$ in $P$, or 
				\item[(b)]  $x,y \in Q$ and $x \leq y$ in $Q$
			\end{enumerate}
			\item The ordinal sum $P \oplus Q$ is defined as the poset on $P \cup Q$, where $x \leq y$ in $P \oplus Q$ if 
			\begin{enumerate}
				\item[(a)] $x,y \in P$ and $x \leq y $ in $P$, or 
				\item[(b)] $x,y \in Q$ and $x \leq y $ in $Q$, or
				\item[(c)] $y \in P$ and $x \in Q$
			\end{enumerate}		
		\end{enumerate}
	\end{definition}
	If the posets $P$ and $Q$ are isomorphic, we write $P \cong Q$.
	
	\section{Tables of Known Values}\label{values}
	\begin{table}[H]
		\centering
		\begin{tabular}{rr}
			$n$ & $\overset{\Delta}{H}(n)$\\ \hline
			1&1 \\
			2&3 \\
			3&14\\
			4&61\\
			5&273\\
			6&1228\\
			7&5631\\
			8&26141\\
			9&123261\\
			10&589251\\
			11&2855815\\
			12&14021493\\
			13&69708557\\
			14&$>$311082529
		\end{tabular}
		\caption{Stationary differences $\overset{\Delta}{H}(n)$ (\seqnum{A376894}).}
		\label{Hdelta}
	\end{table}
	\begin{table}[H]
		\centering
		\begin{tabular}{rr}
			$n$ & $\overset{\infty}{H}(n)$\\ \hline
			0&1 \\
			1&1 \\
			2&4 \\
			3&12\\
			4&47\\
			5&174\\
			6&749\\
			7&3291\\
			8&15675\\
			9&78104\\
			10&411042\\
			11&2261961\\
			12&13009112\\
			13&77860689\\
			14&$>$383637634
		\end{tabular}
		\caption{Limit number of posets $\overset{\infty}{H}(n)$ without isolated points and $n$ arcs in the Hasse diagram (\seqnum{A022016}).}
		\label{Hinf}
	\end{table}
	\begin{landscape}
		\begin{longtable}{r|rrrrrrrrrrrrrr}
			\centering
			\diagbox[height=1cm, width =1cm]{$a$}{$p$} &
			\textbf{1} &
			\textbf{2} &
			\textbf{3} &
			\textbf{4} &
			\textbf{5} &
			\textbf{6} &
			\textbf{7} &
			\textbf{8} &
			\textbf{9} &
			\textbf{10} &
			\textbf{11} &
			\textbf{12} &
			\textbf{13} &
			\textbf{14}\\ \hline
			\endhead
			\textbf{0}  & \textbf{1} & 1 & 1 & 1 & 1  & 1   & 1   & 1    & 1     & 1      & 1       & 1         & 1   & 1       \\ 
			\textbf{1}  & 0 & 1 & 1 & 1 & 1  & 1   & 1   & 1    & 1     & 1      & 1       & 1         & 1   & 1       \\
			\textbf{2}  & 0 & \textbf{0} & \textbf{3} & 4 & 4  & 4   & 4   & 4    & 4     & 4      & 4       & 4         & 4     & 4     \\
			\textbf{3}  & 0 & 0 & 0 & 8 & 11 & 12  & 12  & 12   & 12    & 12     & 12      & 12        & 12  &12       \\
			\textbf{4}  & 0 & 0 & 0 & 2 & \textbf{29} & \textbf{43}  & 46  & 47   & 47    & 47     & 47      & 47        & 47   & 47      \\
			\textbf{5}  & 0 & 0 & 0 & 0 & 12 & 105 & 156 & 170  & 173   & 174    & 174     & 174       & 174   & 174     \\
			\textbf{6}  & 0 & 0 & 0 & 0 & 5  & 92  & 460 & \textbf{670}  & \textbf{731}   & 745    & 748     & 749       & 749      &749  \\
			\textbf{7}  & 0 & 0 & 0 & 0 & 0  & 45  & 582 & 2097 & 2954  & 3212   & 3273    & 3287      & 3290  & 3291  \\
			\textbf{8}  & 0 & 0 & 0 & 0 & 0  & 12  & 487 & 3822 & 10513 & 14196  & \textbf{15323}   & \textbf{15596}     & 15657     & 15671 \\
			\textbf{9}  & 0 & 0 & 0 & 0 & 0  & 3   & 204 & 4514 & 24584 & 55352  & 71548   & 76545     & 77752    & 78025  \\
			\textbf{10} & 0 & 0 & 0 & 0 & 0  & 0   & 71  & 3271 & 40182 & 160398 & 307688  & 381373    & 403979   & 409462  \\
			\textbf{11} & 0 & 0 & 0 & 0 & 0  & 0   & 14  & 1579 & 43365 & 339849 & 1055096 & 1781575   & 2125378   & 2229606 \\
			\textbf{12} & 0 & 0 & 0 & 0 & 0  & 0   & 7   & 561  & 32506 & 513620 & 2792498 & 7044376   & 10723243  & 12368235 \\
			\textbf{13} & 0 & 0 & 0 & 0 & 0  & 0   & 0   & 186  & 17165 & 558180 & 5549235 & 22449776  & 47726348  & 66759843 \\
			\textbf{14} & 0 & 0 & 0 & 0 & 0  & 0   & 0   & 44   & 7361  & 443481 & 8291941 & 56174560  & 178169427  & 328689738\\
			\textbf{15} & 0 & 0 & 0 & 0 & 0  & 0   & 0   & 16   & 2471  & 269258 & 9363668 & 109635845 & 541658175 &1403768258 \\
			\textbf{16} & 0 & 0 & 0 & 0 & 0  & 0   & 0   & 4    & 830   & 129365 & 8163332 & 167196214 & 1326649835 & 5040870856\\
			\textbf{17} & 0 & 0 & 0 & 0 & 0  & 0   & 0   & 0    & 231   & 52171  & 5610656 & 201154432 & 2612415266 & 14997085913\\
			\textbf{18} & 0 & 0 & 0 & 0 & 0  & 0   & 0   & 0    & 73    & 18409  & 3135921 & 193508261 & 4154955648 & 36781063103 \\
			\textbf{19} & 0 & 0 & 0 & 0 & 0  & 0   & 0   & 0    & 18    & 6033   & 1463318 & 151355746 & 5380198148 & 74443114797\\
			\textbf{20} & 0 & 0 & 0 & 0 & 0  & 0   & 0   & 0    & 9     & 1908   & 595863  & 97982560  & 5732689374 & 124919944162 \\
			\textbf{21} & 0 & 0 & 0 & 0 & 0  & 0   & 0   & 0    & 0     & 619    & 217504  & 53635462  & 5085158388 & 175007040480\\
			\textbf{22} & 0 & 0 & 0 & 0 & 0  & 0   & 0   & 0    & 0     & 168    & 75409   & 25410855  & 3804609913 & 206368426296\\
			\textbf{23} & 0 & 0 & 0 & 0 & 0  & 0   & 0   & 0    & 0     & 56     & 24358   & 10709996  & 2435073587 & 206673723878\\
			\textbf{24} & 0 & 0 & 0 & 0 & 0  & 0   & 0   & 0    & 0     & 20     & 8066    & 4118306   & 1354980899 & 177473287015\\
			\textbf{25} & 0 & 0 & 0 & 0 & 0  & 0   & 0   & 0    & 0     & 5      & 2555    & 1485140   & 666949877  & 132025145641\\
			\textbf{26} & 0 & 0 & 0 & 0 & 0  & 0   & 0   & 0    & 0     & 0      & 826     & 513152    & 296092723  & 86035801840\\
			\textbf{27} & 0 & 0 & 0 & 0 & 0  & 0   & 0   & 0    & 0     & 0      & 238     & 172776    & 120852253  & 49712231637\\
			\textbf{28} & 0 & 0 & 0 & 0 & 0  & 0   & 0   & 0    & 0     & 0      & 91      & 56963     & 46293008  &25804712188 \\
			\textbf{29} & 0 & 0 & 0 & 0 & 0  & 0   & 0   & 0    & 0     & 0      & 22      & 18874     & 16899316 & 12205203149  \\
			\textbf{30} & 0 & 0 & 0 & 0 & 0  & 0   & 0   & 0    & 0     & 0      & 11      & 6141      & 5983275  & 5338441123  \\
			\textbf{31} & 0 & 0 & 0 & 0 & 0  & 0   & 0   & 0    & 0     & 0      & 0       & 2018      & 2064635  & 2191738059  \\
			\textbf{32} & 0 & 0 & 0 & 0 & 0  & 0   & 0   & 0    & 0     & 0      & 0       & 623       & 704506   & 856342834  \\
			\textbf{33} & 0 & 0 & 0 & 0 & 0  & 0   & 0   & 0    & 0     & 0      & 0       & 208       & 235107  & 322245931   \\
			\textbf{34} & 0 & 0 & 0 & 0 & 0  & 0   & 0   & 0    & 0     & 0      & 0       & 68        & 78804 &117808997     \\
			\textbf{35} & 0 & 0 & 0 & 0 & 0  & 0   & 0   & 0    & 0     & 0      & 0       & 24        & 25896 & 42116156     \\
			\textbf{36} & 0 & 0 & 0 & 0 & 0  & 0   & 0   & 0    & 0     & 0      & 0       & 6         & 8679   &14779342    \\
			\textbf{37} & 0 & 0 & 0 & 0 & 0  & 0   & 0   & 0    & 0     & 0      & 0       & 0         & 2718    & 5104627   \\
			\textbf{38} & 0 & 0 & 0 & 0 & 0  & 0   & 0   & 0    & 0     & 0      & 0       & 0         & 922    & 1735283    \\
			\textbf{39} & 0 & 0 & 0 & 0 & 0  & 0   & 0   & 0    & 0     & 0      & 0       & 0         & 290     & 583441   \\
			\textbf{40} & 0 & 0 & 0 & 0 & 0  & 0   & 0   & 0    & 0     & 0      & 0       & 0         & 109      & 194263  \\
			\textbf{41} & 0 & 0 & 0 & 0 & 0  & 0   & 0   & 0    & 0     & 0      & 0       & 0         & 26    & 64280     \\
			\textbf{42} & 0 & 0 & 0 & 0 & 0  & 0   & 0   & 0    & 0     & 0      & 0       & 0         & 13      & 21075   \\
			\textbf{43} &0 &0 &0 &0 &0 &0 &0 &0 &0 &0 &0 &0 &0 &6943\\
			\textbf{44} &0 &0 &0 &0 &0 &0 &0 &0 &0 &0 &0 &0 &0 &2225\\
			\textbf{45} &0 &0 &0 &0 &0 &0 &0 &0 &0 &0 &0 &0 &0 &757\\
			\textbf{45} &0 &0 &0 &0 &0 &0 &0 &0 &0 &0 &0 &0 &0 &248\\
			\textbf{46} &0 &0 &0 &0 &0 &0 &0 &0 &0 &0 &0 &0 &0 &80\\
			\textbf{47} &0 &0 &0 &0 &0 &0 &0 &0 &0 &0 &0 &0 &0 &28\\
			\textbf{48} &0 &0 &0 &0 &0 &0 &0 &0 &0 &0 &0 &0 &0 &7\\
			\textbf{49} &0 &0 &0 &0 &0 &0 &0 &0 &0 &0 &0 &0 &0 &0\\
			\caption{Number of posets $H(p,a)$ with $p$ unlabeled points and $a$ arcs in the Hasse diagram (\seqnum{A342447}). The bold numbers indicate where the limit number $\overset{\infty}{H}(n)$ in the difference $H(p,a)- H(p-1,a)$ is reached for the first time (Table~\ref{Hdelta}).}
			\label{normal}
		\end{longtable}
	\end{landscape}

	\begin{landscape}
		\begin{center}
			\begin{longtable}{r|rrrrrrrrrrrrrr}
				
				\diagbox[height=1cm, width =1cm]{$a$}{$p$} &
				\textbf{1} &
				\textbf{2} &
				\textbf{3} &
				\textbf{4} &
				\textbf{5} &
				\textbf{6} &
				\textbf{7} &
				\textbf{8} &
				\textbf{9} &
				\textbf{10} &
				\textbf{11} &
				\textbf{12} &
				\textbf{13} &
				\textbf{14}\\ \hline
				\endhead
				\textbf{0}  & 1 & 0 & 0 & 0 & 0  & 0  & 0   & 0    & 0     & 0      & 0       & 0         & 0     &0     \\
				\textbf{1}  & 0 & 1 & 0 & 0 & 0  & 0  & 0   & 0    & 0     & 0      & 0       & 0         & 0      &0    \\
				\textbf{2}  & 0 & 0 & 3 & 0 & 0  & 0  & 0   & 0    & 0     & 0      & 0       & 0         & 0       &0   \\
				\textbf{3}  & 0 & 0 & 0 & 8 & 0  & 0  & 0   & 0    & 0     & 0      & 0       & 0         & 0     &0     \\
				\textbf{4}  & 0 & 0 & 0 & 2 & 27 & 0  & 0   & 0    & 0     & 0      & 0       & 0         & 0     &0     \\
				\textbf{5}  & 0 & 0 & 0 & 0 & 12 & 91 & 0   & 0    & 0     & 0      & 0       & 0         & 0        &0  \\
				\textbf{6}  & 0 & 0 & 0 & 0 & 5  & 87 & 350 & 0    & 0     & 0      & 0       & 0         & 0     &0     \\
				\textbf{7}  & 0 & 0 & 0 & 0 & 0  & 45 & 532 & 1376 & 0     & 0      & 0       & 0         & 0     &0     \\
				\textbf{8}  & 0 & 0 & 0 & 0 & 0  & 12 & 475 & 3272 & 5743  & 0      & 0       & 0         & 0     & 0     \\
				\textbf{9}  & 0 & 0 & 0 & 0 & 0  & 3  & 201 & 4298 & 19396 & 24635  & 0       & 0         & 0     &0     \\
				\textbf{10} & 0 & 0 & 0 & 0 & 0  & 0  & 71  & 3197 & 36664 & 113734 & 108968  & 0         & 0      &0    \\
				\textbf{11} & 0 & 0 & 0 & 0 & 0  & 0  & 14  & 1565 & 41706 & 292435 & 657864  & 492180    & 0    &0      \\
				\textbf{12} & 0 & 0 & 0 & 0 & 0  & 0  & 7   & 554  & 31931 & 479273 & 2223644 & 3775486   & 2266502   &0 \\
				\textbf{13} & 0 & 0 & 0 & 0 & 0  & 0  & 0   & 186  & 16972 & 540413 & 4952995 & 16241325  & 21518776  & 10598452\\
				\textbf{14} & 0 & 0 & 0 & 0 & 0  & 0  & 0   & 44   & 7317  & 435913 & 7829469 & 47217841  & 114916376 &122051693 \\
				\textbf{15} & 0 & 0 & 0 & 0 & 0  & 0  & 0   & 16   & 2455  & 266743 & 9086436 & 99774543  & 422058803  &792120435\\
				\textbf{16} & 0 & 0 & 0 & 0 & 0  & 0  & 0   & 4    & 826   & 128519 & 8031366 & 158740623 & 1148778277 & 3582725982\\
				\textbf{17} & 0 & 0 & 0 & 0 & 0  & 0  & 0   & 0    & 231   & 51936  & 5557611 & 195406962 & 2402305776 & 12192228986\\
				\textbf{18} & 0 & 0 & 0 & 0 & 0  & 0  & 0   & 0    & 73    & 18336  & 3117269 & 190317665 & 3955460187 &32402651896\\
				\textbf{19} & 0 & 0 & 0 & 0 & 0  & 0  & 0   & 0    & 18    & 6015   & 1457212 & 149873335 & 5225554978 &68853972661 \\
				\textbf{20} & 0 & 0 & 0 & 0 & 0  & 0  & 0   & 0    & 9     & 1899   & 593937  & 97380455  & 5633190375 & 119027215964\\
				\textbf{21} & 0 & 0 & 0 & 0 & 0  & 0  & 0   & 0    & 0     & 619    & 216876  & 53416005  & 5030909731 &169819834744\\
				\textbf{22} & 0 & 0 & 0 & 0 & 0  & 0  & 0   & 0    & 0     & 168    & 75241   & 25334800  & 3778976166 &202508520193 \\
				\textbf{23} & 0 & 0 & 0 & 0 & 0  & 0  & 0   & 0    & 0     & 56     & 24302   & 10685470  & 2424286385 &204212633876\\
				\textbf{24} & 0 & 0 & 0 & 0 & 0  & 0  & 0   & 0    & 0     & 20     & 8046    & 4110184   & 1350837769 &176107385057\\
				\textbf{25} & 0 & 0 & 0 & 0 & 0  & 0  & 0   & 0    & 0     & 5      & 2050    & 1482565   & 665456523 &131354009826 \\
				\textbf{26} & 0 & 0 & 0 & 0 & 0  & 0  & 0   & 0    & 0     & 0      & 826     & 512321    & 295576961 &85738201565 \\
				\textbf{27} & 0 & 0 & 0 & 0 & 0  & 0  & 0   & 0    & 0     & 0      & 238     & 172538    & 120678636  &49590859136\\
				\textbf{28} & 0 & 0 & 0 & 0 & 0  & 0  & 0   & 0    & 0     & 0      & 91      & 56872     & 46235807  &25758244084 \\
				\textbf{29} & 0 & 0 & 0 & 0 & 0  & 0  & 0   & 0    & 0     & 0      & 22      & 18852     & 16880351  &12188246237 \\
				\textbf{30} & 0 & 0 & 0 & 0 & 0  & 0  & 0   & 0    & 0     & 0      & 11      & 6130      & 5977112  &5332438723  \\
				\textbf{31} & 0 & 0 & 0 & 0 & 0  & 0  & 0   & 0    & 0     & 0      & 0       & 2018      & 2062606   & 2189667228 \\
				\textbf{32} & 0 & 0 & 0 & 0 & 0  & 0  & 0   & 0    & 0     & 0      & 0       & 623       & 703883    & 855636277 \\
				\textbf{33} & 0 & 0 & 0 & 0 & 0  & 0  & 0   & 0    & 0     & 0      & 0       & 208       & 234899   & 322010201  \\
				\textbf{34} & 0 & 0 & 0 & 0 & 0  & 0  & 0   & 0    & 0     & 0      & 0       & 68        & 78736   & 117729985   \\
				\textbf{35} & 0 & 0 & 0 & 0 & 0  & 0  & 0   & 0    & 0     & 0      & 0       & 24        & 25872  &42090192    \\
				\textbf{36} & 0 & 0 & 0 & 0 & 0  & 0  & 0   & 0    & 0     & 0      & 0       & 6         & 8673    &14770639   \\
				\textbf{37} & 0 & 0 & 0 & 0 & 0  & 0  & 0   & 0    & 0     & 0      & 0       & 0         & 2718    &5101903   \\
				\textbf{38} & 0 & 0 & 0 & 0 & 0  & 0  & 0   & 0    & 0     & 0      & 0       & 0         & 922    & 1734361    \\
				\textbf{39} & 0 & 0 & 0 & 0 & 0  & 0  & 0   & 0    & 0     & 0      & 0       & 0         & 290   & 583151     \\
				\textbf{40} & 0 & 0 & 0 & 0 & 0  & 0  & 0   & 0    & 0     & 0      & 0       & 0         & 109     & 194154   \\
				\textbf{41} & 0 & 0 & 0 & 0 & 0  & 0  & 0   & 0    & 0     & 0      & 0       & 0         & 26      &64254   \\
				\textbf{42} & 0 & 0 & 0 & 0 & 0  & 0  & 0   & 0    & 0     & 0      & 0       & 0         & 13     &21062    \\
					\textbf{43} &0 &0 &0 &0 &0 &0 &0 &0 &0 &0 &0 &0 &0 &6943\\
				\textbf{44} &0 &0 &0 &0 &0 &0 &0 &0 &0 &0 &0 &0 &0 &2225\\
				\textbf{45} &0 &0 &0 &0 &0 &0 &0 &0 &0 &0 &0 &0 &0 &757\\
				\textbf{45} &0 &0 &0 &0 &0 &0 &0 &0 &0 &0 &0 &0 &0 &248\\
				\textbf{46} &0 &0 &0 &0 &0 &0 &0 &0 &0 &0 &0 &0 &0 &80\\
				\textbf{47} &0 &0 &0 &0 &0 &0 &0 &0 &0 &0 &0 &0 &0 &28\\
				\textbf{48} &0 &0 &0 &0 &0 &0 &0 &0 &0 &0 &0 &0 &0 &7\\
				\textbf{49} &0 &0 &0 &0 &0 &0 &0 &0 &0 &0 &0 &0 &0 &0\\
				\caption{Number of connected posets $H_c(p,a)$ with $p$ unlabeled points and $a$ arcs in the Hasse diagram (\seqnum{A342590}).}\label{conn}
			\end{longtable}
		\end{center}
	\end{landscape}
	\begin{landscape}
		\begin{center}
			\begin{longtable}{r|rrrrrrrrrrrrr}
				
				\diagbox[height=1cm, width =1cm]{$a$}{$p$} &
				\textbf{1} &
				\textbf{2} &
				\textbf{3} &
				\textbf{4} &
				\textbf{5} &
				\textbf{6} &
				\textbf{7} &
				\textbf{8} &
				\textbf{9} &
				\textbf{10} &
				\textbf{11} &
				\textbf{12} &
				\textbf{13} \\
				\textbf{0}  & 1 & 0 & 0 & 0 & 0 & 0  & 0  & 0  & 0   & 0    & 0    & 0     & 0      \\
				\textbf{1}  & 0 & 1 & 0 & 0 & 0 & 0  & 0  & 0  & 0   & 0    & 0    & 0     & 0      \\
				\textbf{2}  & 0 & 0 & 1 & 0 & 0 & 0  & 0  & 0  & 0   & 0    & 0    & 0     & 0      \\
				\textbf{3}  & 0 & 0 & 0 & 2 & 0 & 0  & 0  & 0  & 0   & 0    & 0    & 0     & 0      \\
				\textbf{4}  & 0 & 0 & 0 & 2 & 3 & 0  & 0  & 0  & 0   & 0    & 0    & 0     & 0      \\
				\textbf{5}  & 0 & 0 & 0 & 0 & 2 & 7  & 0  & 0  & 0   & 0    & 0    & 0     & 0      \\
				\textbf{6}  & 0 & 0 & 0 & 0 & 1 & 11 & 10 & 0  & 0   & 0    & 0    & 0     & 0      \\
				\textbf{7}  & 0 & 0 & 0 & 0 & 0 & 5  & 16 & 26 & 0   & 0    & 0    & 0     & 0      \\
				\textbf{8}  & 0 & 0 & 0 & 0 & 0 & 4  & 17 & 68 & 39  & 0    & 0    & 0     & 0      \\
				\textbf{9}  & 0 & 0 & 0 & 0 & 0 & 1  & 9  & 72 & 106 & 107  & 0    & 0     & 0      \\
				\textbf{10} & 0 & 0 & 0 & 0 & 0 & 0  & 5  & 73 & 182 & 386  & 160  & 0     & 0      \\
				\textbf{11} & 0 & 0 & 0 & 0 & 0 & 0  & 0  & 41 & 190 & 675  & 648  & 458   & 0      \\
				\textbf{12} & 0 & 0 & 0 & 0 & 0 & 0  & 1  & 24 & 183 & 1027 & 1600 & 2208  & 702    \\
				\textbf{13} & 0 & 0 & 0 & 0 & 0 & 0  & 0  & 8  & 94  & 1051 & 2547 & 5455  & 3820   \\
				\textbf{14} & 0 & 0 & 0 & 0 & 0 & 0  & 0  & 4  & 63  & 959  & 3585 & 11155 &  12492  \\
				\textbf{15} & 0 & 0 & 0 & 0 & 0 & 0  & 0  & 2  & 23  & 647  & 3522 & 16761 & 27133  \\
				\textbf{16} & 0 & 0 & 0 & 0 & 0 & 0  & 0  & 2  & 16  & 236  & 3404 & 22457 & 51227  \\
				\textbf{17} & 0 & 0 & 0 & 0 & 0 & 0  & 0  & 0  & 5   & 146  & 2419 & 23824 & 73672  \\
				\textbf{18} & 0 & 0 & 0 & 0 & 0 & 0  & 0  & 0  & 3   & 57   & 1821 & 23899 & 98951  \\
				\textbf{19} & 0 & 0 & 0 & 0 & 0 & 0  & 0  & 0  & 0   & 35   & 1018 & 19533 & 107148 \\
				\textbf{20} & 0 & 0 & 0 & 0 & 0 & 0  & 0  & 0  & 1   & 13   & 661  & 15947 & 112153 \\
				\textbf{21} & 0 & 0 & 0 & 0 & 0 & 0  & 0  & 0  & 0   & 8    & 286  & 10517 & 97211  \\
				\textbf{22} & 0 & 0 & 0 & 0 & 0 & 0  & 0  & 0  & 0   & 2    & 181  & 7228  & 83784  \\
				\textbf{23} & 0 & 0 & 0 & 0 & 0 & 0  & 0  & 0  & 0   & 4    & 58   & 3998  &  59871\\
				\textbf{24} & 0 & 0 & 0 & 0 & 0 & 0  & 0  & 0  & 0   & 1    & 46   & 2492  & 44411  \\
				\textbf{25} & 0 & 0 & 0 & 0 & 0 & 0  & 0  & 0  & 0   & 0    & 14   & 1183  & 26479  \\
				\textbf{26} & 0 & 0 & 0 & 0 & 0 & 0  & 0  & 0  & 0   & 0    & 12   & 747   & 17669  \\
				\textbf{27} & 0 & 0 & 0 & 0 & 0 & 0  & 0  & 0  & 0   & 0    & 0    & 324   & 9080   \\
				\textbf{28} & 0 & 0 & 0 & 0 & 0 & 0  & 0  & 0  & 0   & 0    & 3    & 206   &  5877 \\
				\textbf{29} & 0 & 0 & 0 & 0 & 0 & 0  & 0  & 0  & 0   & 0    & 0    & 72    & 2543   \\
				\textbf{30} & 0 & 0 & 0 & 0 & 0 & 0  & 0  & 0  & 0   & 0    & 1    & 58    &  1690  \\
				\textbf{31} & 0 & 0 & 0 & 0 & 0 & 0  & 0  & 0  & 0   & 0    & 0    & 20    & 624    \\
				\textbf{32} & 0 & 0 & 0 & 0 & 0 & 0  & 0  & 0  & 0   & 0    & 0    & 19    & 491    \\
				\textbf{33} & 0 & 0 & 0 & 0 & 0 & 0  & 0  & 0  & 0   & 0    & 0    & 4     & 147    \\
				\textbf{34} & 0 & 0 & 0 & 0 & 0 & 0  & 0  & 0  & 0   & 0    & 0    & 4     & 136    \\
				\textbf{35} & 0 & 0 & 0 & 0 & 0 & 0  & 0  & 0  & 0   & 0    & 0    & 2     & 26     \\
				\textbf{36} & 0 & 0 & 0 & 0 & 0 & 0  & 0  & 0  & 0   & 0    & 0    & 2     & 43     \\
				\textbf{37} & 0 & 0 & 0 & 0 & 0 & 0  & 0  & 0  & 0   & 0    & 0    & 0     & 6      \\
				\textbf{38} & 0 & 0 & 0 & 0 & 0 & 0  & 0  & 0  & 0   & 0    & 0    & 0     & 10     \\
				\textbf{39} & 0 & 0 & 0 & 0 & 0 & 0  & 0  & 0  & 0   & 0    & 0    & 0     & 0      \\
				\textbf{40} & 0 & 0 & 0 & 0 & 0 & 0  & 0  & 0  & 0   & 0    & 0    & 0     & 3      \\
				\textbf{41} & 0 & 0 & 0 & 0 & 0 & 0  & 0  & 0  & 0   & 0    & 0    & 0     & 0      \\
				\textbf{42} & 0 & 0 & 0 & 0 & 0 & 0  & 0  & 0  & 0   & 0    & 0    & 0     & 1     \\
				\caption{Number of connected self-dual posets $H_{\mathrm{sd}}(p,a)$ with $p$ unlabeled points and $a$ arcs in the Hasse diagram (\seqnum{A376633}).}\label{selfdual}
			\end{longtable}
		\end{center}
	\end{landscape}

\end{appendices}

	\bibliography{Handrich_Numbers.bib}
	\bibliographystyle{spmpsci.bst}
\end{document}